   \documentclass[12pt,onecolumn ]{IEEEtran}

 \IEEEoverridecommandlockouts
 
\usepackage{url}
\usepackage{amsmath,amssymb}
\usepackage{amsthm,thmtools}
\usepackage{amsfonts} 

\usepackage[dvips]{epsfig}      
\usepackage{graphicx}
\usepackage{enumerate}
\usepackage{color}                  

\usepackage{verbatim}

\usepackage{amssymb}
\usepackage{amsbsy}
\usepackage{amsmath}
\usepackage{setspace}
\usepackage{url}

 \newcommand{\real}{\operatorname{Re}}
 \newcommand{\imag}{\operatorname{Im}}
 \newcommand{\diag}{\operatorname{diag}}

 \newcommand{\adju}{\operatorname{adj}}

\newcommand*\diff{\mathop{}\!\mathrm{d}}

\declaretheoremstyle[bodyfont=\normalfont]{normalfont}

\declaretheorem[name={Example},qed={\lower-0.3ex\hbox{$\square$}} ] {Example}

\declaretheorem[name={Definition}  ] {Definition}
\declaretheorem[name={Theorem}, style=normalfont] {Theorem}
\declaretheorem[name={Lemma}, style=normalfont ] {Lemma}
\declaretheorem[name={Remark}  ] {Remark}
\declaretheorem[name={Corollary} ,style=normalfont ] {Corollary}
\declaretheorem[name={Assumption} ,style=normalfont ] {Assumption}
\declaretheorem[name={Proposition}, style=normalfont ] {Proposition}

\newcommand {\R}{\mathbb R}

\newcommand {\C}{\mathbb C}

\newcommand{\be}{\begin{equation}}
\newcommand{\ee}{\end{equation}}

\newcommand{\V}{\mathcal V}

\begin{document}
%

\doublespace

\title{On the spectral properties of  nonsingular
matrices that are strictly sign-regular for some order with applications 
to totally positive discrete-time~systems\thanks{The research  of Michael Margaliot is
supported in part 
by  research grants from  the Israel Science Foundation  and the 
 US-Israel Binational Science Foundation.  
}}
 \author{Rola Alseidi,    Michael Margaliot, and  J{\"u}rgen Garloff\thanks{
\IEEEcompsocthanksitem
 Rola Alseidi and J{\"u}rgen Garloff  are with the 
Department of Mathematics and Statistics, University of Konstanz, Germany. 
\IEEEcompsocthanksitem
Michael Margaliot (corresponding author) is  with the School of Electrical  Engineering,
Tel-Aviv University, Tel-Aviv~69978, Israel.
 E-mail: \texttt{michaelm@eng.tau.ac.il}
 }}
 
\maketitle
\begin{center} 
\vspace*{-2cm} 
								Version Date: Oct. 25, 2018
\end{center}
\begin{abstract}
A matrix is called strictly sign-regular of order~$k$ (denoted by~$SSR_k$) if all its~$k\times k$ minors are non-zero and have the same sign. 
For example, totally positive matrices, i.e., matrices with all minors positive, are~$SSR_k$ for all~$k$.
Another important subclass are those that are~$SSR_k$ for all odd~$k$. 
Such matrices have interesting sign variation diminishing properties, and it has been recently shown that
they play an important role in the analysis of certain nonlinear cooperative dynamical systems. 

In this paper, the spectral properties of  nonsingular 
matrices that are~$SSR_k$ for a specific value~$k$ are studied. 
One of the results is 
that the product of the first~$k$ eigenvalues is real and of the same sign as the~$k\times k$ minors,
and that    linear combinations of certain eigenvectors   have specific sign patterns. 
It is then shown how known results for matrices that are~$SSR_k$ for several values of~$k$ can be derived from these spectral properties.

Using these theoretical results,  the notion of a totally positive discrete-time  system~(TPDTS) is introduced.
This may be regarded as the discrete-time  analogue of the important notion of a totally positive differential system, introduced by Schwarz in~1970. 
The asymptotic behavior of  time-invariant and time-varying~TPDTSs is analyzed, and it
 is shown that every trajectory of a  periodic time-varying~TPDTS  
converges to a periodic solution. 
\end{abstract}
\begin{center}
\end{center}

\section{Introduction}

A matrix (not necessarily square) is called   \emph{sign-regular of order~$k$} (denoted by~$SR_k$) if all its minors of order~$k$ are non-negative or all are non-positive. 
It is called  \emph{strictly sign-regular of order~$k$}  (denoted by~$SSR_k$)
if it is   sign-regular of order~$k$, and all the minors of order~$k$ are non-zero. In other words, all minors of order~$k$ are 
non-zero and have the same sign. A matrix is called 
\emph{sign-regular}~($SR$)
  if it is~$SR_k$ for all~$k$, and
\emph{strictly sign-regular}~($SSR$)
  if it is~$SSR_k$   for all~$k$.
 	
The most prominent examples of	$SR$ [$SSR$]
	  are totally nonnegative \emph{TN} [totally positive \emph{TP}]
 matrices, that is,
matrices  with all minors nonnegative [positive].
Such matrices have applications in a number of fields including 
approximation theory, economics, probability theory,  computer aided geometric design and more~\cite{total_book,gk_book,pinkus}.  

A very important property  of~$SSR$ matrices is  that multiplying a vector by such a matrix can only decrease the number of sign variations. 
To explain this variation diminishing property~(VDP), we introduce some notation.
We use small letters to denote column vectors. If $y \in \mathbb{R}^n$ is such a vector then~$y'$ denotes its transpose. 
For ~$y \in \mathbb{R}^n$, we use
 $s^-(y)$ to denote  the number of sign variations
	in~$y$ after deleting all its zero entries,
	and~$s^+(y)$ to denote  the maximal possible number of sign variations
	in~$y$ after each zero entry is replaced by either~$+1$ or~$-1$. 
	For example, for~$n=4$ and~$y=\begin{bmatrix} 1& -1 & 0  &  -\pi  \end{bmatrix}'$,
	$s^-(y)=1$ and~$s^+(y)=3$. 
	Obviously,
	\[
	0\leq s^-(y) \leq s^+(y)\leq n-1 \text{ for all } y\in\R^n.
	\]

The first important results on the VDP of matrices were obtained by Fekete \cite{fekete1912} and Schoenberg \cite{Schoenberg1930}. Later on, Grantmacher and Krein \cite[Chapter V]{gk_book} elaborated rather completely the various forms of VDPs and worked out the spectral properties of $SR$ matrices. Two examples of such VDPs are: if~$A\in\mathbb{R}^{n\times m}$ ($m \leq n$) 
 is~$SR$ and of rank $m$ then
\begin{equation*}
s^-(Ax)\leq s^-(x) \text{ for all } x \in \mathbb{R}^m,
\end{equation*}
whereas if $A$ is   $SSR$ then  
\begin{equation*}
s^+(Ax)\leq s^-(x) \text{ for all } x \in \mathbb{R}^m \setminus \{0\}.
\end{equation*} 
 The more recent literature on $SR$ matrices focuses on the recognition and factorization of matrices with special signs of their minors, see, e.g., \cite{huang2013test} and the references therein.

There is a renewed interest in such VDPs in the context of dynamical 
 systems. Indeed, 
  Reference~\cite{fulltppaper}   shows
 that 
powerful results on the asymptotic behavior of the solutions of
 continuous-time
nonlinear cooperative and tridiagonal dynamical systems 
can be derived using the fact that the transition matrix of a certain linear time-varying system (called the variational system)
is~TP for all time~$t$. In other words, the linear system is a    totally positive differential system~(TPDS)~\cite{schwarz1970}.
These transition matrices are real, square, and nonsingular. 
In a recent paper~\cite{CTPDS}, it is shown that the   transition matrix satisfies a~VDP 
with respect to the \emph{cyclic} number of sign variations if and only if 
it is~$SSR_k$ for all odd~$k$ and all~$t$.

Some of the spectral properties of $SR$ matrices can be extended to matrices which are $SR_k$ for all $k$ up to a
 certain order~\cite [Chapter V]{gk_book}. In this paper, we study the spectral properties of nonsingular matrices which are $SSR_k$ for a specific value of of $k$. Such matrices are only rarely considered in the literature, see, e.g., \cite{Pena}, where a test for an $n\times k$ matrix with $k <n$ to be $SSR_k$ is presented.  
Let~$\epsilon_k  \in \{-1,1\}$ denote the sign  of the minors of order $k$, with convention $\epsilon_0:=1$. 
We prove that the strict sign-regularity of order~$k$ implies
 that the product of the first~$k$ eigenvalues is real and has the sign~$\epsilon_k$,
and that  
certain  eigenvectors satisfy a special sign pattern.   
Then we show how to extend these  results to obtain the spectral properties of
 matrices that 
are~$SSR_k$ for several values of~$k$, for example, for 
 all odd~$k$.

 The theoretical  results are used to derive a
new class of dynamical systems 
called \emph{totally positive discrete-time  systems}~(TPDTSs). 
This is the discrete-time 
 analogue of the important notion of a~TPDS.
We analyze the asymptotic properties of TPDTSs and show that they 
entrain to periodic excitations. This result may be regarded as the analogue of an important result of Smith~\cite{periodic_tridi_smith}
on entrainment in a continuous-time periodic nonlinear cooperative systems with a special Jacobian. In particular, this result
 implies that any trajectory of a time-invariant~TPDS either escapes every compact set or converges to an equilibrium.

The main novel contributions of this paper include: (1)~the spectral analysis of matrices that are strictly-sign regular for some order~$k$;
(2)~the introduction of a new class of discrete-time linear time-varying dynamical systems called~TPDTSs; and (3)~the analysis of the asymptotic behavior of discrete-time \emph{nonlinear}
 time-varying dynamical systems whose variational equation is~TPDTS. In particular, we prove   that such nonlinear systems entrain to a periodic excitation.

 The remainder of this paper is organized as follows.
The next section provides  known definitions and reviews results that will be used later on.
In Section~\ref{sec:main}, we present our main results and  in Section \ref{sec:aplic}  we
apply these results to 
introduce and analyze~TPDTSs. The final section concludes and outlines possible  topics  for further research.

\section{Preliminaries}
In this section, we
 collect   several known definitions and results that will be used later on.

Given a square matrix~$A\in\R^{n\times n}$ and~$p\in\{1,\dots,n\}$, consider the
$\binom{n}{p}^2$
 minors of~$A$ of order $p$. 
Each minor is defined by a set of~$p$ row indexes
\be\label{eq:rind}
1\leq i_1<i_2<\dots<i_p\leq n
\ee
 and~$p$ column indexes
\be\label{eq:cind}
1\leq j_1<j_2<\dots<j_p\leq n.
\ee
This minor 
is denoted by~$A(\alpha|\beta)$, where~$\alpha:=\{i_1,\dots,i_p\}$ and~$\beta:=\{j_1,\dots,j_p\}$ (with a mild abuse of notation, we will  regard 
these sequences as sets). We suppress the curly brackets if we enumerate the indexes explicitly.  

The next result, known as Jacobi's identity, provides 
 information  on the relation between the minors of
a nonsingular matrix~$A\in\R^{n\times n}$ 
and the minors of~$A^{-1}$. For a sequence~$\alpha=\{i_1,\dots,i_p\}$, let~$\bar \alpha$
denote the sequence~$\{1,\dots,n\}\setminus \alpha$ which we consider as ordered increasingly.

\begin{Proposition}\label{prop:jacobi}\cite[Section 1.2]{total_book}
Let~$A\in \R^{n\times n }$ be a nonsingular matrix, and put~$B:=A^{-1}$. 
Pick~$p\in\{1,\dots,n\}$. 
Then for any two sequences~$\alpha=\{ i_1, i_2,\dots, i_p \} $
and~$\beta=\{ j_1, j_2,\dots, j_p \} $ satisfying~\eqref{eq:rind} and~\eqref{eq:cind}
we have
\[
								B(\alpha,\beta)=(-1)^s \frac{A( \bar \beta,\bar \alpha  ) }{\det(A)}, 
\] 
with~$s:=i_1+\dots+i_p+j_1+\dots+j_p$. 
\end{Proposition}

 Let~$D_{\pm 1}\in\R^{n\times n}$ denote the diagonal matrix with diagonal 
entries~$(1,-1,1,-1,\dots,(-1)^{n+1})$, and let~$\adju(A)$ denote the adjugate of~$A$.
Then~$D_{\pm 1}\adju(A)D_{\pm 1}^{-1}$
is the unsigned adjugate  of~$A$. Proposition~\ref{prop:jacobi} yields the following result. 
\begin{Corollary}\label{prop:corojac}
Let~$A\in \R^{n\times n }$ be a nonsingular matrix.
Suppose that for some~$k\in\{1,\dots,n-1\}$
all the minors of order~$k$ of~$A$ are non-zero and have the same sign. 
Then all minors of order~$n-k $ of the unsigned adjugate of~$A$ are non-zero and have the same sign. 
\end{Corollary}

The $SSR_k$ property is closely related to VDPs.
The following well-known result demonstrates this.
\begin{Proposition}\label{prop:sum1}\cite[Chapter V, Theorem 1]{gk_book}
	Consider a   matrix~$U\in \R^{n\times m }$  with~$n>m $. 
The following  two conditions are equivalent:
\begin{enumerate}[(i)]
\item	\label{cond:cisold} For any~$x\in \R^m\setminus\{ 0 \}$,  we have 
\be\label{eq:suimold}
					s^+(Ux)\leq m-1.
\ee
\item			\label{cond:minoold}
The matrix~$U$ is $SSR_m$, that is, all minors of the form 
\be\label{eq:minors}
U( i_1\;\dots\; i_m| 1\; \dots \;m), \text{ with }1\leq i_1<i_2<\dots<i_m\leq n,
\ee
 are non-zero and have the same sign. 
\end{enumerate}
\end{Proposition}

 Note that the assumption that~$n>m$ cannot be weakened.
For example, if we take~$n=m=2$ 
and a square matrix~$U\in\R^{2\times 2} $
then condition~\eqref{cond:cisold} obviously holds, yet 
 condition~\eqref{cond:minoold}  only holds if~$U$ is nonsingular, so the conditions are not equivalent in this case.

 It was recently shown that
for square  
matrices 
  the~$SSR_k$ property is equivalent to a non-standard~VDP. 
\begin{Theorem}\label{thm:ssrpp1}~\cite[Theorem 1]{CTPDS}
Let~$A\in\R^{n\times n}$ be a nonsingular matrix. 
  Pick~$k \in \{1,\dots,n\}$. 
Then the following  two conditions are equivalent:
\begin{enumerate}[(i)]
\item	\label{cond:onedip} For any vector~$x\in\R^n\setminus\{0\}$ with~$ s^-(x) \leq k-1$, we have 
\[
					s^+(Ax)\leq k-1.
\]
\item			\label{cond:secdip}
$
A 
$
is~$SSR_{k}$. 
\end{enumerate}
\end{Theorem}

 We emphasize that condition~\eqref{cond:onedip} in 
Theorem~\ref{thm:ssrpp1} does not  assert  that~$s^-(x)\leq k-1$
implies that~$s^+(Ax)\leq s^-(x)$,
but only that~$s^+(Ax)\leq k-1$.  

\section{Main Results}\label{sec:main}

We consider from here on a nonsingular matrix~$A\in\R^{n\times n}$. 
  We  order its eigenvalues such that
\be\label{eq:orlam}
						|\lambda_1|\geq |\lambda_2| \geq \dots \geq |\lambda_n|>0,
\ee
with complex conjugate  eigenvalues appearing  in consecutive pairs
(we   say, with a mild abuse of notation,
 that~$z \in \C$ is \emph{complex}  if~$z\not = \bar z$, where~$\bar z$ denotes
 the complex conjugate of~$z$). 
Let
\be\label{eq:vis}
v^1,\; v^2 , \dots, v^n\in\C^n  
\ee	
	 denote the eigenvectors corresponding to the~$\lambda_i$'s.
	We always assume that every~$v^i$ is not purely imaginary.
	Indeed, otherwise  we can replace~$v^i$ 
	by~$\imag(v^i)$ that is a real eigenvector. 
	The fact that~$A$ is real implies  that if~$v^i$ is complex then its real and imaginary parts can be chosen as linearly independent. 
	
	 Define a set of real vectors~$u^1,u^2,\dots,u^n \in \R^n$ by going through the~$v^i$'s as follows.
 If~$v^1$ is real then we put~$u^1:=v^1$ and  proceed to examine~$v^2$.
 If~$v^1$ is complex (and whence~$v^2=\bar v^1$) then we put~$u^1:=\real(v^1)$, $u^2:=\imag(v^1)$ and proceed to examine~$v^3$, and so on. Let~$U:=\begin{bmatrix} u^1&\dots&u^n \end{bmatrix} \in\R^{n\times n }$. 

Suppose that for some~$i,k$ the eigenvector~$v^i$ is real and
$v^k$ is complex. Then 
  is not difficult to show that
since~$A$ is real and nonsingular,
the real vectors~$v^i,\real(v^k),\imag(v^k)$ are linearly independent.

	 Note that if~$v^i,v^{i+1}\in\C^n$ is a complex conjugate pair and~$c\in \C\setminus\{0\}$ is complex then 
\[
			c   v^i+\bar c v^{i+1} =2( \real(c) \real(v^i) -\imag(c)\imag(v^i)) \in \R^n\setminus\{0\},
\]
 and  by choosing an appropriate complex~$c \in \C\setminus\{0\}$ we can get any nonzero real linear combination of the real vectors~$ \real(v^i)$ and~$\imag(v^i) $.

We say that a set~$c_p,\dots,c_k\in \C$, $p\leq k$, \emph{matches} the set $v^p,\dots,v^k$ of consecutive eigenvectors~\eqref{eq:vis}  if the~$c_i$'s are not all zero and for every~$i$ if  the vector~$v^i$ is real then~$c_i$ is real,
	and if~$v^i,v^{i+1}$ is a complex conjugate pair then~$c_{i+1}=\bar c_i$.  
	In particular, this implies that~$\sum_{i=p}^k   c_i v^i \in  \R^n. $

In order to prove our main results, we need the following auxiliary result that is a generalization of 
Proposition~\ref{prop:sum1} to the 
case of eigenvectors of a real and nonsingular matrix. 
We use the notation~$j$ for the imaginary unit ($j^2=-1$). 
\begin{Proposition}\label{prop:gensum1} 
Consider the set of~$n$ vectors~$v^1,\dots,v^n \in \C^n$.
Define~$V\in \C^{n\times n }$ by
\[
V:=\begin{bmatrix} v^1& v^2&\dots&v^n \end{bmatrix}.
\]
The following  two conditions are equivalent:
\begin{enumerate}[(i)]
\item	\label{cond:cis} For any~$c_1,\dots,c_m\in \C$  that 
match~$v^1,\dots,v^m$, with~$m\leq n$, we have
\be 
					s^+(\sum_{i=1}^m c_i v^i)\leq m-1 .
\ee
\item			\label{cond:mino}
Let~$q \in \C^{ \binom{n}{m}} $ be the vector that 
 contains all  the minors of order~$m$ 
of~$V$ of the form $V(i_1 \dots i_m|1 \dots m)$, with $1 \leq i_1<i_2<\dots<i_m \leq n$, 
arranged in the lexicographic order.  
Then there exists an integer~$k$ such that 
the vector~$j ^k q$ is real
and   all its entries are positive.  
\end{enumerate}
\end{Proposition}
{\sl Proof.} 
Consider first the case~$m=1$.
In this case, the definition of a matching set means that~$v^1$
is real. 
Condition~\eqref{cond:cis} becomes
\be\label{eq:cond1tr}
s^+(v^1)=0
\ee
and condition~\eqref{cond:mino} becomes
\[
v^1=j^k q
\]
for some integer~$k\geq 0$ and~$q$ a real vector with all entries positive. Since~$v^1$ is real, this means that all the entries of~$v^1$ are either all positive or all negative. Thus, in this case the two conditions are indeed equivalent.  

Suppose  now that~$m=2$.
Then two cases are possible.

\noindent {\sl Case 1.}
Both~$v^1$ and~$v^2$ are real. 
Then by the definition of a matching set, \eqref{cond:cis} 
becomes: for any~$d_1,d_2\in\R$, that are not both zero,
\[ 
					s^+(d_1 v^1 +d^2 v^2)\leq 1 .
\]
On the other-hand, \eqref{cond:mino} becomes: the vector~$q \in \C^{ \binom{n}{2}} $  that 
 contains all  the minors of order $2$ of $V$
 of the form $V(i_1 \; i_2|1 \; 2)$, with $1 \leq i_1<i_2  \leq n$, arranged in the lexicographic order,  
is real and   all its entries are non-zero and have   the same sign.  Now the equivalence of the two conditions follows from Proposition~\ref{prop:sum1}. 

\noindent {\sl Case 2.}
Both~$v^1$ and~$v^2$ are complex. In this case,~$v^2=\real(v^1)-j\imag(v^1)$,  
and~\eqref{cond:cis} 
becomes: for any~$d_1,d_2\in\R$, that are not both zero,
\[ 
					s^+(d_1 \real(v^1) +d_2 \imag(v^1))\leq 1 .
\]
On the other-hand, \eqref{cond:mino} becomes: let~$q \in \C^{ \binom{n}{2}} $  be the vector that 
 contains all  the minors of order $2$ of the matrix 
\[
V=\begin{bmatrix}  \real(v^1) +j\imag(v^1) & \real(v^1) -j\imag(v^1)   \end{bmatrix} 
\]
of the form $V(i_1 \; i_2|1 \;2)$, with $1 \leq i_1<i_2  \leq n$, arranged in the lexicographic order.  
Then  there exists an integer~$k$ such that the vector~$j ^k q$ is real,
and   all its entries are positive.
  
Let~$U:= \begin{bmatrix}  \real(v^1)   &  \imag(v^1)   \end{bmatrix}$. Then  
\[
V(i_1 \; i_2|1\; 2)=(-2j)  U(i_1 \; i_2|1 \;2)
\]
for all~$1 \leq i_1<i_2  \leq n$.
Since~$U$ is a real matrix,~\eqref{cond:mino} becomes: let~$q \in \R^{ \binom{n}{2}} $  be the vector that 
 contains all  the minors of order $2$ of~$U$ 
of the form $U(i_1 \; i_2|1 \;2)$, with $1 \leq i_1<i_2  \leq n$, arranged in the lexicographic order.  
Then~$ q$ is real,
and   all its entries are non-zero and have   the same sign.
 Now the equivalence of the two conditions follows from Proposition~\ref{prop:sum1}. 
For any~$m>2$ the proof follows similarly by decomposing the~$v^i$s into sets of real and complex conjugate pairs and then using the results in the cases~$m=1$ and~$m=2$ described above.~\hfill{$\square$}

\subsection{Matrices that are $SSR_k$ for some value $k$}	
We now state our first  main result that describes the spectral properties
of a nonsingular matrix that is~$SSR_k$ for some value~$k$.
\begin{Theorem}\label{thm:styr}
Suppose that~$A\in\R^{n\times n}$ is nonsingular and~$SSR_k$ for
 some value~$k$, with~$k \in \{1,\dots,n-1\}$.  
Then the following properties hold:
\begin{enumerate}[(i)]
\item The product~$\lambda_1 \lambda_2 \dots \lambda_k$ is real, and 
\be\label{eq:porfisp}
\epsilon_k  \lambda_1 \lambda_2 \dots \lambda_k>0.
\ee
\item  The eigenvalues satisfy the inequality
\be\label{eq:lamkkp1}
|\lambda_k|>|\lambda_{k+1}|.
\ee
\item 
Pick~$1\leq p\leq k$,~$k+1\leq q \leq n$, and~$c_p,\dots,c_q  \in \C$ 
  such that
$c_p,\dots,c_k$   [$c_{k+1},\dots,c_q$] match   
   the eigenvectors~$v^p,\dots,v^k$ [$v^{k+1},\dots,v^q$] of~$A$. Then  
\be\label{eq:suim}
					s^+(\sum_{i=p}^k c_i v^i)\leq k-1 ,
\ee
and
\be\label{eq:lpoyt}
					s^-(\sum_{i=k+1}^q c_i v^i)\geq k .
\ee
\item\label{item:lind}
Let~$u^1,\dots,u^{n}$ be the set of real vectors
constructed from~$v^1,\dots,v^n$ as described above. 
Then~$u^1,\dots,u^{k}$ 
are linearly independent. In particular, if~$v^1,\dots,v^k$
are real then they are linearly independent.
\end{enumerate}
\end{Theorem}
\begin{Remark}
  Roughly speaking, equations~\eqref{eq:suim}
	and~\eqref{eq:lpoyt} imply 
		that the first~$k$ [last~$n-k$] 
	eigenvectors of~$A$ have a sign pattern with a ``small'' [``large''] number of sign changes.
In particular, 
 for any~$i\leq k$ we have that if~$v^i \in \R^n$   then~$ s^+(  v^i)\leq k-1$,
and if~$v^{i-1},v^i$ is a complex conjugate pair then for any~$d_1,d_2\in \R$, that are not both zeros,
 we have~$s^+( d_1 \real (v^i ) + d_2\imag(v^i))\leq k-1$. 

Similarly,
for any~$j\geq k+1$ we have that if~$v^j \in \R^n$   then~$ s^-(  v^j)\geq k$,
and if~$v^{j},v^{j+1}$ is a complex conjugate pair then for any~$d_1,d_2\in\R$, that are not both zero,
 we have~$s^-( d_1 \real (v^j ) + d_2\imag(v^j))\geq k$.
\end{Remark}

{\sl Proof of Theorem~\ref{thm:styr}.}
Let~$r:=\binom{n}{k}$.  Recall that the~$k$th \emph{multiplicative  compound matrix}~$A^{(k)}$ is the~$r \times   r$ matrix
that 
includes all the minors of order~$k$ of~$A$ ordered lexicographically. By Kronecker's theorem, see, e.g.,  \cite[p. 132]{pinkus},
the eigenvalues~$\zeta_i   $, $i=1,\dots,r$,
 of~$A^{(k)}$ are all the~$k$ products of~$k$ eigenvalues of~$A$, that is,
\begin{align}\label{eqn:zetai}
										\zeta_1  &= \lambda_1\lambda_2\dots \lambda_{k-1} \lambda_k,\nonumber \\
										\zeta_2   &= \lambda_1\lambda_2\dots  \lambda_{k-1} \lambda_{k+1},\nonumber \\
										&\vdots\\
										\zeta_r   &= \lambda_{n-k+1} \lambda_{n-k+2} \dots  \lambda_{n-1} \lambda_{n}.\nonumber 
\end{align}
Combining this with~\eqref{eq:orlam}, \eqref{eqn:zetai}  implies that
\be\label{eq:pouuit}
	|\zeta_1   | \geq  |	\zeta_2   | \geq\dots\geq 	 |\zeta_r   | >0.
\ee
Since~$A$ is~$SSR_k$, all the entries in~$A^{(k)}$ are non-zero and have the  sign~$\epsilon_k$.
Thus, all the entries of~$\epsilon_k A^{(k)}$ are positive. 
Perron's theorem implies that 
$\epsilon_k \zeta_1 $ is real and positive with a corresponding entry-wise positive eigenvector~$w$ (that is unique up to scaling), and that 
\[
	 \epsilon_k \zeta_1   >  |	\zeta_2    | \geq\dots\geq 	 |\zeta_r    | >0.
\]
Using~\eqref{eqn:zetai}, we conclude that~$\epsilon_k \lambda_1\lambda_2\dots\lambda_k$ is real and positive and that
\[
|\lambda_k |>|\lambda_{k+1}|.
\]

To prove the sign patterns of the eigenvectors, let~$V:=\begin{bmatrix} v^1&\dots&v^n \end{bmatrix}$ and~$D:=\diag(\lambda_1,\dots,\lambda_n)$,
so that~$AV=VD$.
  Define $q\in\C^r$ by $q_\alpha:=V(\alpha|1,2,\dots,k)$, where $\alpha$ is running over all $k$-tuples ~$1\leq i_1<\dots<i_k\leq n$  and the components of~$q$ are ordered lexicographically.
Then the definitions of~$ A^{(k)}$ and~$q$ imply that 
\begin{align*}
(A^{(k)}q)_\alpha &= \sum_{\beta} A(\alpha,\beta)V(\beta|1,2,\dots,k),
\end{align*}
where the summation is over all~$k$-tuples~$
\beta$. 
Let~$B:=AV$. Applying  the Cauchy-Binet formula for the minors
 of the product of two matrices, e.g., \cite[Theorem 1.1.1]{total_book}, yields 
$(A^{(k)}q)_\alpha = B(\alpha|1,2,\dots,k)$. Since~$B=VD$, a further application of the Cauchy-Binet formula results in  
\begin{align*}
(A^{(k)}q)_\alpha& =    \sum_{\beta} V(\alpha,\beta)D(\beta|1,2,\dots,k)   \\
&= V(\alpha| 1,2,\dots,k) \lambda_1\lambda_2\dots\lambda_k\\
&= \zeta_1  q_\alpha, 
\end{align*}
 where the second equation follows from the fact that~$D$ is a diagonal matrix.
Since this holds for any entry of the vector~$A^{(k)}q$, 
we conclude that~$ q$ is an eigenvector of~$\epsilon_k A^{(k)}$ corresponding to its Perron root~$\epsilon \zeta_1$.
Thus,  there exists~$\eta \in \C\setminus\{0\}$ such that~$q=\eta  w$, where~$w\in\R^r$ is an entry-wise positive vector.
Using the fact that the complex vectors in~$V$ appear in   conjugate pairs, 
and arguing as in the proof of Proposition~\ref{prop:gensum1}, we have that~$\eta=j^k$ for some integer~$k$. Now application of Proposition~\ref{prop:gensum1} yields  that
for  any~$c_1,\dots,c_k\in \C$  that match~$v^1,\dots,v^k$ we have
\[					s^+(\sum_{i=1}^k c_i v^i)\leq k-1 
\]
which proves~\eqref{eq:suim}.

To prove~\eqref{eq:lpoyt},
note that~$Av^i=\lambda_i v^i$ and the identity~$\adju(A) A= \det(A)I$ yield
\[
				 	(	D_{\pm 1}\adju(A) D_{\pm 1}^{-1} )D_{\pm 1} \lambda_i v^i= \det(A) D_{\pm 1}v^i. 
\]
In other words, the eigenvalues~$\eta_i$ of the unsigned adjugate of~$A$,
ordered such that~$|\eta_1|\geq|\eta_2|\geq\dots\geq |\eta_n|>0$, 
are
\[
\eta_1=		 \frac{\det(A)}{\lambda_n} ,\;
\eta_2= \frac{\det(A)}{\lambda_{n-1}} ,\; \dots,	\;	    \eta_n=\frac{\det(A)}{\lambda_1} ,
\]
with corresponding eigenvectors
\[
 z^1:=D_{\pm 1}v^n, \; z^2:= D_{\pm 1}v^{n-1},\; \dots ,z^{n}:= D_{\pm 1}v^1.
\]
Since~$A$ is nonsingular and~$SSR_k$, Corollary~\ref{prop:corojac}
implies that
 all minors of order~$n-k $ of the unsigned adjugate of~$A$ are non-zero and have the same sign, that is, $	D_{\pm 1}\adju(A) D_{\pm 1}^{-1} $ is~$SSR_{n-k}$. 
By Equation~\eqref{eq:suim} (that has already been proved), this means that
for any~$c_1,\dots,c_{n-k}\in\C$ that match~$z^1,\dots,z^{n-k}$, we have~$s^+(\sum_{i=1}^{n-k} c_i z^i )\leq n-k-1 $, i.e.,
$
s^+(D_{\pm 1 } \sum_{i=k+1}^n c_{n-i+1}  v^i)\leq n-k-1 .
$
Combining this with  the identity
\[
							s^+(D_{\pm 1}x)+s^-(  x)=n-1      \text{ for all } x\in\R^n,
\]
	see, e.g.,~\cite[p. 88]{total_book}, yields~\eqref{eq:lpoyt}.

To prove~\eqref{item:lind} note that since~$k\leq n-1$,
\eqref{eq:suim} yields
$	s^+(\sum_{i=1}^k c_i v^i)\leq n-2$. This means that~$	 \sum_{i=1}^k c_i v^i  \not =0$ for any matching~$c_1,\dots,c_k\in\C$, that is,~$	 \sum_{i=1}^k d_i u^i \not =0$
for any~$d_1,\dots,d_k \in \R$ that are not all zero.					
 This completes the proof of  Theorem~\ref{thm:styr}.~\hfill{$\square$}

\begin{Example}\label{exa:pott}
Consider the 
matrix
\be\label{eq:ant}
A=\begin{bmatrix}  1& 2 & 0 &0 \\
											0&1&1&0 \\ 0 & 0&2 &1\\1&0&0 &2
\end{bmatrix}.
\ee
 It is straightforward to verify that this matrix is nonsingular, and that
 all minors of order $3$ are  positive,  
so~$A$ is~$SSR_3$ with~$\epsilon_3=1$
(but not~$SSR_1$ nor~$SSR_2$). 
 Its eigenvalues   are (all numerical values in this paper 
are to four-digit accuracy)
\[
\lambda_1=2.7900,\;
\lambda_2=1.5000 + 1.0790 j, \;
\lambda_3=1.5000 - 1.0790 j,\;
\lambda_4= 0.2100,
	\]
and thus~$\lambda_1\lambda_2\lambda_3$ is real and positive, and~$|\lambda_3|>|\lambda_4|$.
The matrix of corresponding eigenvectors is
\begin{align*}
V&=\begin{bmatrix} v^1 & v^2 & v^3 &v^4 \end{bmatrix}\\
&=\begin{bmatrix}
							 0.79&		 -0.5 + 1.079 j&	-0.5 - 1.079 j&	-1.79\\
							0.7071& -0.7071 & -0.7071 & 0.7071\\
							1.266& -0.3536 - 0.763 j &-0.3536 + 0.763 j& -0.5586\\
							1&1&1&1
\end{bmatrix},
\end{align*}
and thus
\begin{align*}
									U&:=\begin{bmatrix} v^1& \real(v^2) &\imag(v^2) & v^4 \end{bmatrix}\\
									&=\begin{bmatrix}
							 0.79  &		 -0.5 &	 1.079 &	-1.79\\
							0.7071 & -0.7071 & 0 & 0.7071\\
							1.266  & -0.3536  &- 0.763 & -0.5586\\
							1&1&0&1
\end{bmatrix}. 
\end{align*}

Calculating the  vector~$q$ that contains all 
  minors in the form~$V(\alpha|1,2,3)$ yields~$q=-jw=j^3 w$, with
\[
w:= \begin{bmatrix}
		1.7049 &  3.0518 &  5.4629   &2.1580 
\end{bmatrix}  '. \] 

A calculation of all the minors of the form~$U(\alpha|1,2,3)$ gives the values
$
\begin{bmatrix}
		0.8524& 1.5259& 2.7315&  1.0790
\end{bmatrix}  '. $
Since these are all positive, application of Proposition~\ref{prop:sum1}
to the submatrix containing  the first three columns of~$U$ gives that
 for any~$d_1,d_2,d_3\in\R$, that are not all zero,
\[
s^+(  d_1 v^1 +d_2 \real(v^2)+ d_3 \imag(v^2))\leq 2,
\]
which agrees with~\eqref{eq:suim}. Furthermore, it holds that~$s^-(v^4)\geq 3$
 which agrees with~\eqref{eq:lpoyt}.
\end{Example}

So far, we have considered matrices that are~$SSR_k$ for a single value of~$k$.
Our next goal is to demonstrate 
 that Theorem~\ref{thm:styr} can be used as a basic building block in the analysis
of matrices 
 that are~$SSR_k$ for \emph{several} values of~$k$, by simply applying
 Theorem~\ref{thm:styr} for every such~$k$.

\subsection{Matrices that are~$SSR_k$ for two consecutive values of~$k$}
The next result analyzes matrices that are~$SSR_k$   for two
consecutive values of~$k$.
\begin{Corollary}\label{coro:consek}
Suppose that~$A\in\R^{n\times n}$ is nonsingular,~$SSR_i$ and~$SSR_{i+1}$
for
 some value~$i$, with~$i \in \{1,\dots,n-2\}$. 
Then the following properties hold:
\begin{enumerate}[(i)]
\item \label{item:poiit}
The signed eigenvalue $\epsilon_i \epsilon_{i+1} \lambda_{i+1}$ is real and 
positive.
\item  The eigenvalues satisfy the inequalities
\be\label{eq:lamkk12}
|\lambda_i|>|\lambda_{i+1}|>|\lambda_{i+2}|.
\ee
\item The eigenvector~$v^{i+1}$ can be chosen as a real vector and 
\be\label{eq:poddd}
s^-(v^{i+1})=s^+(v^{i+1})=i.
\ee
 Furthermore, for any~$p,q$ 
with~$1\leq p\leq i$,~$i+2\leq q \leq n$, and~$c_p,\dots,c_i$ 
[$c_{i+2},\dots,c_q$] $\in \C$  that  match   
   the eigenvectors~$v^p,\dots,v^i$ [$v^{i+2},\dots,v^q$] we have   
\be\label{eq:portyy}
					s^+(\sum_{\ell =p}^i c_\ell  v^\ell )\leq i-1 ,\; s^+(\sum_{\ell =p}^{i+1} c_\ell  v^\ell )\leq i ,
\ee
and
\be\label{eq:piytr}
					s^-(\sum_{\ell =i+1}^q c_\ell  v^\ell )\geq i ,\; s^-(\sum_{\ell =i+2}^{q} c_\ell  v^\ell )\geq i +1.
\ee
\item The vectors~$u^1,\dots,u^{i+1}$ are linearly independent. 
\end{enumerate}
\end{Corollary}
{\sl Proof.}
 By Theorem~\ref{thm:styr},~$\epsilon_i \lambda_1 \lambda_2\dots\lambda_i$ 
and~$\epsilon_{i+1} \lambda_1 \lambda_2 \dots \lambda_{i+1}$ are real and positive which yields~(\ref{item:poiit}). Thus,~$v^{i+1}$ can be chosen as a real vector.
 Using~\eqref{eq:lamkkp1}
with~$k=i$ and with~$k=i+1$ gives~(\ref{eq:lamkk12}).
Inequalities~\eqref{eq:suim} and~\eqref{eq:lpoyt} imply~\eqref{eq:portyy}
and~\eqref{eq:piytr},  and this implies
\[
i\leq s^-( v^{i+1}) \leq s^+( v^{i+1})\leq i.
\]
The last statement of the corollary  follows immediately from Theorem~\ref{thm:styr}
and this completes the proof.~\hfill{$\square$}

\subsection{Matrices that are~$SSR$ }
Known results on the spectral structure of~TP matrices, see, e.g,~\cite[Chapter 5]{total_book}, \cite[Chapter 5]{pinkus}, (and, more generally, of $SSR$ matrices \cite[Chapter V]{gk_book}),
 follow immediately from Corollary~\ref{coro:consek}. Indeed, suppose that~$A\in\R^{n\times n}$ is~$SSR$. Then 
it is~$SSR_i$ and~$SSR_{i+1}$ for all~$i\in\{1,\dots,n-1\}$ and thus
Corollary~\ref{coro:consek}
implies that  the eigenvalues~$\lambda_i$, $i=1,\dots, n-1$, are real and that
\be\label{eq:eibgty}
|\lambda_1|>|\lambda_2|>\dots>|\lambda_n|>0.
\ee
Since~$\det(A) =\lambda_1 \dots \lambda_n$ is real, this implies that~$\lambda_n$ is also real. 
Pick indexes~$p,q$ with~$1\leq p \leq q\leq n$, and~$c_p,c_{p+1},\dots,c_q\in \R$
 such that~$v:=\sum_{i=p}^q c_i v^i  \not =0 $. Then 
since~\eqref{eq:portyy}
and~\eqref{eq:piytr}
	 hold  for any~$i$, we conclude that  
\[
   p-1\leq s^-(v) \leq  s^+(v) \leq q-1.
\]
In particular, taking~$p=q$ yields
\[
   s^-(v^p) = s^+(v^p) =p-1\text{ for all } p\in\{1,\dots,n\}.
\]

\subsection{Matrices that are~$SSR_k$ for all odd~$k$}
We now analyze the spectral properties of matrices that are~$SSR_k$ for all odd~$k$. 
To explain why such matrices are interesting, we  first review  their VDPs. 
For a vector~$y\in\R^n$, let 
	\begin{align}\label{eq:defsc-}
	s^-_c(y)& :=\max_{ i\in\{1,\dots,n\} } 
	s^-(\begin{bmatrix}  y_i & \dots &y_n &y_1& \dots&y_i
	\end{bmatrix}').
	\end{align}
	Intuitively speaking, this corresponds to placing
	the entries of~$y$ along a circular ring so that~$y_n$ is followed by~$y_1$, 
	then counting~$s^-$ starting  from 
	any entry along the ring, and finding the maximal value. 
	
	For example, for~$y=\begin{bmatrix}  0&-1&0&2&0&3\end{bmatrix}'$,
	$s^-_c(y)=s^-(\begin{bmatrix}   -1&0&2&0&3&0&-1 \end{bmatrix}') = 2$. 
	Similarly, let 
	\begin{align*}
	 s^+_c(y) & :=\max_{ i\in\{1,\dots,n\} }  
				s^+(\begin{bmatrix}  y_i & \dots &y_n &y_1& \dots&y_i
	\end{bmatrix}'),
	\end{align*}
	but here if~$y_i=0$ then in the calculation of~$s^+(\begin{bmatrix}  y_i & \dots &y_n &y_1& \dots&y_i
	\end{bmatrix}')$ \emph{both}~$y_i$'s are replaced by either~$1$ or~$-1$. 
	In our above example, we have $s^+_c(y)=4$. 
	Clearly,~$s^-_c(y) \leq s^+_c(y)$ for all~$y\in\R^n$,
	and~$s^-_c(y),s^+_c(y)$ are invariant 
	under  cyclic shifts of the vector~$y$.

	There is a simple and useful
	relation between the cyclic and non-cyclic     number of sign variations of a vector, namely,
for any vector~$x$,
\[
				s^-_c(x)=\begin{cases} 
							s^-(x), & \text{ if }s^-(x) \text{ is even},\\
							s^-(x)+1, & \text{ if }s^-(x) \text{ is odd}, 
							\end{cases}
\]
and, similarly, 
\[
				s^+_c(x)=\begin{cases} 
							s^+(x), & \text{ if } s^+(x) \text{ is even},\\
							s^+(x)+1, & \text{ if } s^+(x) \text{ is odd},
									\end{cases}
\]
see, e.g.,~\cite[Chapter 5]{karlin_tp}, where also other useful results of the cyclic variations of sign can be found.
 This implies in particular that for any vector~$x$,
\be\label{eq:nodeeevem}
s^-_c(x) ,s^+_c(x) \in \begin{cases} \{0,2,4,\dots,n \}, & \text{ if } n \text{ is even},\\
																		\{0,2,4,\dots,n-1 \}, & \text{ if } n \text{ is odd}.
\end{cases} 
\ee

It was shown in~\cite{CTPDS} that a nonsingular matrix~$A\in\R^{n\times n}$
satisfies the \emph{cyclic~VDP}, i.e.,
\be\label{eq:mainc1}
					s^+_c(Ax)\leq s^-_c(x)  \text{ for all } x\in\R^n\setminus\{0\},
\ee
 if and only if~$A$ is $SSR_k$ for all \emph{odd}~$k$ in the range~$\{1, \dots,n\}$.

The proof of the next result uses Theorem~\ref{thm:styr} 
to analyze the spectral properties of such matrices. 
\begin{Theorem}\label{thm:ytre}
	Suppose that~$A\in \mathbb{R}^{n\times n}$
	is nonsingular and~$SSR_k$ for all odd~$k$ in the range $\{1, \dots,n\}$.  
	Then 
	\begin{enumerate}[(i)]
		\item \label{item:prt} The eigenvalue   $\lambda_1$ of $A$  is  simple, real, with~$\epsilon_1 \lambda_1>0$.
		\item \label{item:rrrt} The algebraic multiplicity of any eigenvalue of $A$ is not greater than~$2$, and
  the eigenvalues satisfy the inequalities
		$|\lambda_1|>|\lambda_2|\geq|\lambda_3|>|\lambda_4|\geq|\lambda_5|>\dots$ .
		\item \label{item:bnprt} For every even~$k$ in the range $ \{ 2, \dots,n-1\}$, the inequality $\epsilon_{k-1}\epsilon_{k+1}\lambda_{k}\lambda_{k+1}	>0$ holds.
	 
		\item \label{item:prtyu} If $n$ is even then $\lambda_n$ is real,   and~$\epsilon_{n-1}\det(A)\lambda_n>0$.
		\item \label{item:eigvec} For any~$i$, the eigenvectors have the following properties:
		if~$v^{2i+1}$ is real then~$s^+(v^{2i+1})\leq 2i$, 
		and if~$v^{2i},v^{2i+1}$ is a complex conjugate pair then for any matching~$c_1,c_2\in \C$
			\[
			s^+( c_1 v^{2i}+c_2v^{2i+1} )\leq 2i. 
			\]
		Also, if~$v^{2i}$ is real then~$s^-(v^{2i})\geq 2i-1$, 
		and if~$v^{2i},v^{2i+1}$ is a complex conjugate pair then for any matching~$c_1,c_2\in \C$
			\[
			s^-( c_1 v^{2i}+c_2v^{2i+1} )\geq 2i-1. 
			\]
		Furthermore, the vectors~$u^1,\dots,u^p$, with~$p$ the largest odd number satisfying~$p\leq n$,  
	  are linearly independent. 
	\end{enumerate}
\end{Theorem}
{\sl Proof}.
	 Statement (\ref{item:prt}) follows from the fact that~$A$ is~$SSR_1$
	and by Perron's theorem. 
By~\eqref{eq:lamkkp1} it follows that 
	\[
	|\lambda_1|>|\lambda_2|, \; |\lambda_3|>|\lambda_4|,\; |\lambda_5|>|\lambda_6|\;\dots
	\]
	which implies~(\ref{item:rrrt}). 
	By Theorem~\ref{thm:styr}, the products
	\[
	\epsilon_1 \lambda_1,\;\epsilon_3 \lambda_1\lambda_2\lambda_3,\;\epsilon_5\lambda_1\lambda_2\lambda_3 \lambda_4\lambda_5, \;\dots
	\]
	are all real and positive, which implies  $\epsilon_{k-1}\epsilon_{k+1}\lambda_{k}\lambda_{k+1}	>0$ and thus (\ref{item:bnprt}). 
	To prove~(\ref{item:prtyu}), note that if~$n$ is even then~$n-1$ is odd, 
	so~$A$ is~$SSR_{n-1}$. This implies that~$\beta:=\epsilon_{n-1}
	\lambda_1\lambda_2\dots\lambda_{n-1}$ is real and positive, and from the fact that~$\epsilon_{n-1}\det(A)=\beta \lambda_n$
	the claim follows.
	The results on the sign pattern of the eigenvectors follow from~\eqref{eq:suim}
	and~\eqref{eq:lpoyt}, and the linear independence of~$u^1,\dots,u^p$
	follows from   Theorem~\ref{thm:styr} (\ref{item:lind}).~\hfill{$\square$}

	\begin{Example}
	Consider a nonsingular matrix~$A\in\R^{3\times 3}$ with all entries positive.
	Then~$A$ is~$SSR_1$, with~$\epsilon_1=1$, and~$SSR_3$. 
	Theorem~\ref{thm:ytre} implies that $\lambda_{1}$ is positive,
 and  $s^-(v^1)=s^+(v^1)=0$, and also that only one of the following two cases is possible. 
	
	\noindent {\sl Case 1.} 
	 The eigenvalues $\lambda_2, \lambda_3$ are both real, with
	\[
					\lambda_1>|\lambda_2 |\geq |\lambda_3|
	\]
	
	and~$s^-(v^2)\geq 1$.   The eigenvectors~$v^1,v^2,v^3$ are linearly independent.

	\noindent {\sl Case 2.} 
 The eigenvalues $\lambda_{2}, \lambda_{3}$ are complex with $\lambda_2=\bar \lambda_3$. Then
	\[
					\lambda_1>|\lambda_2 |= |\lambda_3|
	\]
 and~$s^-(d_1 \real(v^2)+d_2\imag(v^2))\geq 1$ for all~$d_1,d_2\in \R$  that are not voth zero.
The vectors~$v^1,\real(v^2),\imag(v^2)$ are linearly independent. 
	\end{Example}

	
	\subsection{Matrices that are $SSR_k$ for all even $k$}
We start with a result on the spectral properties of matrices which are $SSR_k$ for all even $k$. Its proof is parallel to one of Theorem \ref{thm:ytre}.
\begin{Corollary}
	
	Let $A\in \mathbb{R}^{n,n}$ be nonsingular and $SSR_k$ for all even $k$ in the range $\{2, \dots,n\}$. Then

	\begin{enumerate}[(i)]
		\item \label{item:rrrt1} The algebraic multiplicity of any eigenvalue of $A$ is not greater than~$2$, and
		the eigenvalues satisfy the inequalities
		$|\lambda_1|\geq|\lambda_2|>|\lambda_3|\geq|\lambda_4|>|\lambda_5|>\dots$ .
		\item \label{item:bnprt1} For every odd~$k$ in the range $ \{ 1,\dots,n-1\}$, the inequality   $\epsilon_{k-1}\epsilon_{k+1}\lambda_{k}\lambda_{k+1}	>0$ holds. 
		
		\item \label{item:prtyu1} If $n$ is odd then $\lambda_n$ is real,   and~$\epsilon_{n-1}\det(A)\lambda_n>0$.
		\item \label{item:eigvec1} For any~$i$, the eigenvectors have the following properties:
		if~$v^{2i}$ is real then~$s^+(v^{2i})\leq 2i-1$, 
		and if~$v^{2i-1},v^{2i}$ is a complex conjugate pair then for any matching~$c_1,c_2 \in \mathbb{C}$
		\[
		s^+( c_1 v^{2i-1}+c_2v^{2i} )\leq 2i-1. 
		\]
		Also, if~$v^{2i-1}$ is real then~$s^-(v^{2i-1})\geq 2i-2$, 
		and if~$v^{2i},v^{2i+1}$ is a complex conjugate pair then for any matching~$c_1,c_2 \in \mathbb{C}$
		\[
		s^-( c_1 v^{2i-1}+c_2v^{2i} )\geq 2i-2. 
		\]
					Furthermore, the vectors~$u^1,\dots,u^p$, with~$p$ the largest even number satisfying~$p\leq n$,  
	  are linearly independent. 
			\end{enumerate}
\end{Corollary}
	We now apply Theorem~\ref{thm:ssrpp1} to show that the matrices which are~$SSR_k$ for all even $k$ also possess a certain VDP.
\begin{Definition}
	Let $x \in \R^n$. Define
		$s^+_o(x):$=$\left\{\begin{matrix}
		s^+(x), \;&\mbox{if }s^+(x)\;\mbox{is odd},\\ 
		s^+(x)+1,\;&\mbox{if } s^+(x)\;\mbox{is even.}\end{matrix}\right.$
	\end{Definition}
	 Note that~$s^+_o(x)$ is the smallest odd number greater than or equal to~$s^+(x)$.
We define analogously~$s^-_o(x)$ by replacing the 
superscript~$+$ by~$-$. 
\begin{Definition}
A matrix~$A \in \mathbb{R}^{m\times n}$ is said to have  the \emph{odd VDP}  if 
\begin{align}\label{eq:ovdp}
s^+_o(Ax) \leq s^-_o(x) \text{ for all }  x \in \mathbb{R}^n \setminus  \{0\}.
\end{align}
  \end{Definition}

\begin{Theorem}\label{thm:pouut}
	Let  $A \in \mathbb{R}^{n\times n}$ be a nonsingular matrix. The following two statements are equivalent: 
		\begin{enumerate}[(i)]
\item $A$ has the odd VDP. \label{con 1 odd}
	\item \label{con 2 odd} The matrix~$A$ is~$SSR_k$ 
	for all even~$k$ in the range $\{2, \dots,n\}$.
\end{enumerate}
\end{Theorem}

{\sl Proof.}
Assume that condition \eqref{con 1 odd} holds. 
Pick an \emph{odd} number $k \in\{1,\dots,n-1\}$ and pick~$x \in \mathbb{R}^n \setminus
 \{0\}$ such that $s^-_o(x) \leq k$. Then~$s^-(x) \leq k$ . Condition (\ref{con 1 odd}) yields $s^+_o(Ax)\leq k$. Theorem~\ref{thm:ssrpp1} implies that~$A$ is~$SSR_{k+1}$. Since~$k$ is an arbitrary odd number in~$\{1,\dots,n-1\}$, we conclude that condition (\ref{con 2 odd}) holds. 

To prove the converse implication, suppose that condition \eqref{con 2 odd} holds. Pick $x \in \mathbb{R}^n \backslash \{0\}.$ Let $p$ be such that $s_o^-(x)=2p-1$ and thus $s^-(x)\leq 2p-1$. If $2p-1=n$ then clearly~\eqref{eq:ovdp} holds; so we may assume that $2p-1 \leq n-1$. Condition (\ref{con 2 odd}) implies in particular that $A$ is $SSR_{2p}$ and application of Theorem  \ref{thm:ssrpp1} yields $s^+(Ax) \leq 2p-1$, hence $s^+_o(Ax) \leq 2p-1=s^-_o(x).$~\hfill{$\square$}

The next section describes an application
of the   properties analyzed above to discrete-time dynamical systems. 
We begin with linear systems whose transition matrix is~TP, and then show how these can be applied  to analyze nonlinear systems. 
	
\section{Applications to Discrete-Time Dynamical Systems}\label{sec:aplic}

\subsection{Linear systems}
Consider the   linear   time-varying   system:
\be\label{eq:dlti}
														x(i+1)=A(i)x(i), \quad x(0)=x_0,
\ee
with~$A(i)\in \R^{n\times n}$ for all~$i\geq 0$, and~$x_0\in\R^n$.

Our goal is to develop for  the discrete-time system~\eqref{eq:dlti} an analogue of the important notion
 of a~TPDS, derived by Schwarz~\cite{schwarz1970} for continuous-time systems. The main idea is to require 
every~$A(i)$ to satisfy a VDP. 

\begin{Lemma}\label{lem:dltiser}
							Suppose that there exists~$k\in\{1,\dots,n-1\}$ such that for every~$i\geq 0$
							the
							matrix~$A(i)$ is nonsingular and~$SSR_k$.
							Pick~$x_0 \in \R^n\setminus\{0\}$ such that~$s^-(x_0) \leq k-1$.
							Then the solution of~\eqref{eq:dlti} satisfies
\[
							 s^+(x(i))\leq k-1 \text{ for all } i\geq 1.
\]
\end{Lemma}
{\sl Proof. } Fix~$i\geq 1$. Then~$x(i)= A(i-1) \dots A(1) A(0)x_0$.
The matrix~$ A(i-1) \dots A(1) A(0)$ is nonsingular and~$SSR_k$, as it is 
 the product of nonsingular and~$SSR_k$ matrices. 
Theorem~\ref{thm:ssrpp1}  implies that~$s^+(x(i)) \leq k-1$.~\hfill{$\square$}

Let~$\V:=\{y\in\R^n: s^-(y)=s^+(y)\}$. For example, for~$n=3$ the vector~$y=\begin{bmatrix} 1& 0 & -1
\end{bmatrix}'\in \V$, as~$s^-(y)=s^+(y)=2$.
It is not difficult to show that
\[
\V=\{y\in \R^n: y_1\not =0, y_n \not =0, \text{ and if }  y_i =0 \text { for some }
i\in\{2,\dots,n-1\} \text{ then } y_{i-1}y_{i+1}<0 \}.
\]		
 
\begin{Theorem}\label{thm:mainds} 
		Suppose that for every~$i\geq 0$ 
		the matrix~$A(i)$ is~$SSR$.
		Then 
		\begin{enumerate}[(i)]
									\item			For any~$x_0\in \R^n\setminus\{0\}$ the
solution of~\eqref{eq:dlti} satisfies
\be\label{eq:tryuue}
							s^+(x(i+1))\leq s^-(x(i)) \text{ for all } i\geq 0,
\ee
and~$x(i) \in \V$ for all~$i\geq 0 $ except perhaps for up to~$n-1$ values of~$i$. \label{item:exce}
\item \label{item:thewoo} 
Pick two different initial conditions~$x_0,\bar x_0\in\R^n $,
and let~$x(q),\bar x(q)$ denote the corresponding solutions 
of~\eqref{eq:dlti} at time~$q$.
Then there exists~$m\geq 0$ such that
\be\label{eq:poyyr}
x_1(q)\not =\bar x_1(q) \text{ for all }q\geq m.
\ee 
\end{enumerate}
\end{Theorem}
{\sl Proof.} Pick~$x_0\in \R^n\setminus\{0\}$. 
If~$s^-(x_0)=n-1$ then clearly~\eqref{eq:tryuue}  holds, 
so we may assume that~$s^-(x_0)<n-1$.
Let~$q\in\{1,\dots,n-1\}$ be  
such that~$s^-(x_0)=q-1$. 
Since~$A(0)$ is nonsingular and~$SSR_q$, 
Lemma~\ref{lem:dltiser}
implies that~$s^-(x(1))\leq s^+(x(1))\leq q-1=s^-(x_0)$. 
Let~$p\leq q$ be such that~$s^-(x(1))=p-1$. Since~$A(1)$ is nonsingular and~$SSR_p$, 
 Lemma~\ref{lem:dltiser}
implies that~$s^-(x(2))\leq s^+(x(2))\leq p-1=s^-(x(1))$. 
Proceeding in this fashion  proves~\eqref{eq:tryuue}.

The analysis above shows that the mappings~$i\to s^+(x(i))$
and~$i\to s^-(x(i))$ are nonincreasing, and if~$x(i)\not \in \V$, that is,~$s^-(x(i))<s^+(x(i))$ then~$s^+(x(i+1))<s^+(x(i))$.   
Since~$s^+$ takes values in~$\{0,\dots,n-1\}$, this  proves~\eqref{item:exce}.

To prove~\eqref{item:thewoo}, let~$z(i):=x(i)-\bar x(i)$.
Then~$z(0)\not =0$,~$z(i+1)=A(i)z(i)$, and~\eqref{item:exce} implies that there exists~$m\geq 0$ such that~$z(p)\in\V$ for 
all~$p\geq m$. In particular,~$z_1(p)\not =0$ for all~$p\geq m$,
 and this completes the proof.~\hfill{$\square$}

\begin{Remark}
Suppose that~$A(i)=I$ for all~$i\geq 0$. 
Pick~$x_0 \in\R^n \setminus\{0\}$ such that~$x_0\not \in \V$.
Then  
the solution of~\eqref{eq:dlti}   satisfies~$x(i) \not \in \V$ for all~$i\geq 0$.
Thus, 
 Theorem~\ref{thm:mainds} does not hold if we weaken
 the hypothesis of the theorem 
to:~``$A(i)$ is nonsingular 
and~$SR$ for all~$i\geq 0$''. 
\end{Remark}

For the analysis of periodic discrete-time nonlinear systems below, we need to strengthen the result in~\eqref{eq:poyyr}. 
To do this, we use a  result from~\cite{gk_book}.
Recall that two vectors~$v,w\in\R^n$ are said to \emph{oscillate in the same way} if~$s^-(v)=s^-(w)$ and the  
  first non-zero entry in~$v$ and in~$w$ have  the same sign. 
	For example,~$v=\begin{bmatrix} 0&0&1&0&-2\end{bmatrix}'$
	and~$w=\begin{bmatrix} 5&0&-2&-3&-3\end{bmatrix}'$ oscillate in the same way. 
	Theorem~5 in~\cite[p.~254]{gk_book} implies that~$A\in \R^{n\times n}$ is~TN if and only if for any~$x\in\R^n$ the vector~$y:=Ax$ 
	satisfies~$s^-(y)\leq s^-(x)$ and if~$s^-(y)= s^-(x)$ then $x$ 
	and~$y$ oscillate in the same way. 

\begin{Definition}
We call~\eqref{eq:dlti} a \emph{totally positive discrete-time system}~(TPDTS)
  if~$A(i)$ is~TP for all~$i\geq 0$. 
\end{Definition}	
This is the analogue of the notion of a~TPDS
 defined by Schwarz~\cite{schwarz1970}
 for continuous-time systems.  Note that if~\eqref{eq:dlti} is~TPDTS then, in particular, every 
entry of~$A(i)$ is positive for all~$i$, so the system is cooperative~\cite{hlsmith}.

\begin{Theorem}\label{thm:stronger} 
		Suppose that~\eqref{eq:dlti}  is~TPDTS.
		Then  
		\begin{enumerate}[(i)]
									\item	\label{item:porr}		For any~$x_0\in \R^n\setminus\{0\}$ the
solution of~\eqref{eq:dlti} satisfies
\be\label{eq:tryuuext}
							s^+(x(i+1))\leq s^-(x(i)) \text{ for all } i\geq 0,
\ee
and~$x(i) \in \V$ for all~$i\geq 0 $ except perhaps for up to~$n-1$ values of~$i$. \label{item:exceqq}
\item \label{item:thewooty} 
Pick two different initial conditions~$x_0,\bar x_0\in\R^n  $,
and let~$x(q),\bar x(q)$ denote the corresponding solutions of~\eqref{eq:dlti} at time~$q$.
Then there exists~$r\geq 0$ such that either
\[
x_1(q)>\bar x_1(q) \text{ for all }q\geq r 
\]
or
\[
x_1(q)<\bar x_1(q) \text{ for all }q\geq r.
\]
\end{enumerate}
\end{Theorem}
{\sl Proof.} Since a~TP matrix is~SSR, statement~\eqref{item:porr}   follows immediately from
Theorem~\ref{thm:mainds}. To prove~\eqref{item:thewooty}, 
recall that~$z(i):=x(i)-\bar x(i)$ satisfies~$z(i)=A(i)z(i)$,
 and that there exists~$m\geq 0$ such that~$z_1(p)\not = 0 $ for all~$p\geq m$. By the analysis above, there 
   exists~$\ell\geq 0$ such that
\[
s^-(z(i))=s^+(z(i)) =s^-(z(i+1)) =s^+(z(i+1)) \text{ for all } i\geq \ell. 
\]
In particular, we obtain
\be\label{eq:possw}
s^-(z(i))=s^-(A(i)z(i))\text{ for all } i\geq \ell.
\ee
Let~$r:=\max(m,\ell)$. Pick~$i\geq r$. 
Then    the first non-zero entry in~$z(i)$ and in~$z(i+1)=A(i)z(i)$ is the first entry. 
Since~$A(i)$ is~TP (and thus~TN), \eqref{eq:possw}
 implies that~$z_1(i)$ and~$z_1(i+1)$ have the same sign. Since~$i \geq r$ is arbitrary, this proves~\eqref{item:thewooty}.~\hfill{$\square$}

 \begin{Remark}
Consider the special case, where~\eqref{eq:dlti}
 is time-invariant, that is,~$A(i)=A$ for all~$i\geq 0$, with~$A$ a~TP matrix.
In this case, it is possible to give a simpler proof for the ``eventual monotonicity'' result in~\eqref{item:thewooty}. Indeed, 
let~$\lambda_1>\dots>\lambda_n>0$ be the eigenvalues of~$A$ with corresponding
 eigenvectors~$v^i\in\R^n$, $i=1,\dots,n$. 
Write~$z_0=x_0-\bar x_0$ as
$
z_0 =\sum_{\ell=1}^n c_\ell v^\ell .
$
Since~$z_0\not =0$, there exists a minimal index~$p$ such that~$c_p \not =0$, that is, 
$z_0 =\sum_{\ell=p}^n c_\ell v^\ell $ and thus
\be\label{eq:pcse}
z(i)=A^i z_0=\sum_{\ell=p}^n c_\ell (\lambda_\ell) ^i v^\ell.
\ee
Recall that~$s^-(v^p)=s^+(v^p)$, so in particular the first entry of~$v^p$ is not zero.
Since~$\lambda_p>\lambda_{p+1}>\dots>\lambda_n$, it follows from~\eqref{eq:pcse} that for any~$i$ sufficiently  
large, the first entry of~$z(i)$
 is not zero and has the same sign as
the first entry of~$c_pv^p$.
\end{Remark}
 
\subsection{Nonlinear systems}
Now consider the time-varying nonlinear discrete-time dynamical system
\be\label{eq:nlol}
x(i+1)= f(i,x(i)),
\ee
with~$f$ continuously differentiable with respect to~$x$.
Let~$J(i,x) :=\frac{\partial}{\partial x} f(i,x) $ 
denote its Jacobian with respect to~$x$. 
We assume throughout 
that
 the  trajectories of~\eqref{eq:nlol} evolve
 on a convex and compact set~$\Omega \subset \R^n$.

Pick two different initial conditions~$x_0,\bar x_0\in\Omega$
and let~$z(i):=x(i)-\bar x(i)$. 
Then
\begin{align}\label{eq:zdyonp}
z(i+1)&=  f(i,x(i))- f(i,\bar x(i))\nonumber \\
&= \int_0^ 1  \frac{d}{dr} f (i,r x(i)+(1-r) \bar x(i))  \diff r\nonumber \\
&=  \int_0^ 1 \frac{\partial }{\partial x} 
 f(i,r x(i)+(1-r) \bar x(i))  \frac{\partial}{\partial r} 
 (r x(i)+(1-r) \bar x(i)) \diff r  \nonumber \\
&=\left ( \int_0^ 1  J(i,r x(i)+(1-r) \bar x(i)) \diff r \right ) 
 z(i).
\end{align}
This is a discrete-time variational system, as it describes how a  variation 
 in the initial condition propagates with time.

The next assumption guarantees that~\eqref{eq:zdyonp} is~TPDTS.
\begin{Assumption}\label{assum:trp}
						The matrix~$ \int_0^ 1  J(i,r a+(1-r) b ) \diff r $ is~TP for any~$a,b \in \Omega $ and any~$i\geq0$.
\end{Assumption} 

To illustrate an application of Theorem~\ref{thm:stronger}, we assume that there exists a~$T>0$ such that 
the nonlinear system~\eqref{eq:nlol} 
is~$T$-periodic. 
\begin{Assumption}\label{assum:TPER}
						The vector field in~\eqref{eq:nlol}
						satisfies~$ f(i+T,a) =f(i,a)$ for all~$a\in \Omega$ and all~$i\geq 0$. 
\end{Assumption}

We  now state the main result in this section.
We say that a solution~$x(i)$ of~\eqref{eq:nlol} 
is a \emph{$T$-periodic solution} if~$x(i+T)=x(i)$ for all~$i\geq 0$. 

\begin{Theorem}\label{thm:nonltprt}
If Assumptions~\ref{assum:trp} and~\ref{assum:TPER}
 hold then for any~$a\in \Omega$ the
 solution~$x(i,a)$ of~\eqref{eq:nlol}
 converges to a $T$-periodic solution as $i\to\infty $. 
\end{Theorem}

If we view the~$T$-periodic vector field in~\eqref{eq:nlol}
as representing a 
periodic excitation   then 
Theorem~\ref{thm:nonltprt}
implies that the system \emph{entrains} 
to the excitation, as every solution converges to a periodic solution with the same period as the excitation. 
Entrainment is important in many natural and artificial systems. Proper  functioning of various organisms requires internal processes to entrain to the 24h solar day. 
Synchronous generators must entrain to the frequency  of the grid. For more details, see, e.g.,~\cite{2017arXiv171007321M,RFM_entrain,entrain2011}.  

{\sl Proof of Theorem~\ref{thm:nonltprt}.}
Pick~$a \in \Omega$.
Let~$x(i,a)$ denote the solution of~\eqref{eq:nlol}
 at time $i$ with~$x_0=a$. If~$x(i,a)$ is~$T$-periodic then there is nothing to prove. Thus, we may assume that~$y(i):=x(i,a)$
and~$w(i):=x(i+T,a)$ are not identical and the $T$-periodicity of the vector field implies that both~$y,w$ are solutions of the dynamical system.
 
 Theorem~\ref{thm:stronger}  
 implies that there exists~$m\geq 0$ such that either~$w_1(i)>y_1(i)$ or~$w_1(i)<y_1(i)$ for all~$i\geq m$. 
Without loss of generality, we may assume  that the first of these two inequalities holds and this yields
\be\label{eq:pdocf}
x_1((k+1)T,a) > x_1(kT,a) \text{ for all } k\geq m.
\ee

Define
\[
\omega_T(a):=\{ v\in \Omega:  \lim_{k\to\infty} x(n_kT,a)=v \text{ for some sequence } n_1<n_2<\dots    \}.
\]
Since the solutions remain in the compact set~$\Omega$,
the set~$\omega_T(a)$ is not empty. 
If~$v \in \omega_T(a)$   then~$x(T,v)=v$, so the solution emanating from~$v$ is~$T$-periodic. Thus, to complete the proof  we need to show that~$\omega_T(a)$ is a singleton.   
We assume on the contrary that   there exist~$p,q\in\Omega$, with~$p\not =q$, such that~$p,q \in  \omega_T(a)$.
Then there exist sequences~$\{n_k\}_{k=1}^\infty$ and~$\{m_k\}_{k=1}^\infty$ such that
\[
p= \lim_{k\to \infty } x(n_k T,a) \text{ and }  q= \lim_{k\to \infty } x(m_k T,a) . 
\]
Passing to subsequences, we may assume that~$n_k<m_k<n_{k+1}$
for all~$k$. 
Now~\eqref{eq:pdocf} yields~$x_1(n_kT,a) <x_1(m_k T ,a )<x_1(n_{k+1}T,a)$ for all~$k$ sufficiently large and taking~$k\to\infty$ yields~$p_1=q_1$. We conclude that
any two points in~$\omega_T(a)$ have the same first coordinate.

Consider the trajectories emanating from~$p$ and~$q$, that is,
$x(i,p)$ and~$x(i,q)$. Since~$p,q\in \omega_T(a)$ and~$\omega_T(a)$ is an invariant set,
$x(kT,p),x(kT,q)\in\omega_T(a)$ for all~$k$, so
\[
			x_1(kT,p)=x_1(kT,q) \text{ for all } k.
\] 
However, since~$p\not =q$  this contradicts the eventual  monotonicity property described in Theorem~\ref{thm:stronger}. 
This contradiction proves that~$\omega_T(a)$ is a singleton which completes the proof.~\hfill{$\square$}

\begin{figure*}[t]
 \begin{center}
  \includegraphics[scale=0.6]{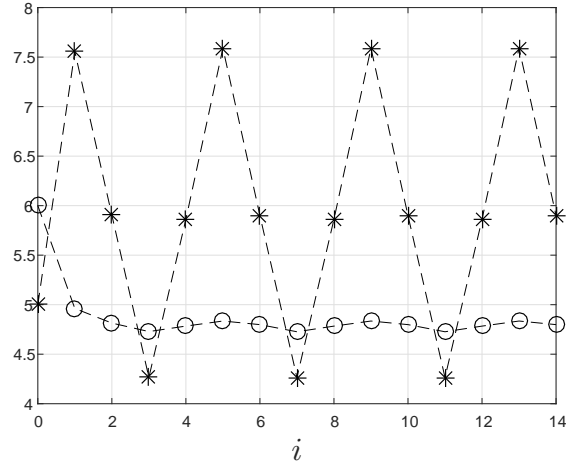}
\caption{$x_1(i)$ (marked by asterisks) 
 and $x_2(i)$ (marked by circles) as a function of~$i$ for the system in Example~\ref{exa:pore}. }\label{fig:imp}
\end{center}
\end{figure*}

 \begin{Example}\label{exa:pore}
Consider the nonlinear system
\begin{align}\label{eq:exalp}
x_1(i+1)&=\frac{c_{11}(i) x_1(i)}{1+x_1(i)}+ \frac{c_{12}(i) x_2(i)}{1+x_2(i)},\nonumber \\
x_2(i+1)&=\frac{c_{21}(i) x_1(i)}{1+x_1(i)}+ \frac{c_{22}(i) x_2(i)}{1+x_2(i)}.
\end{align}
Its Jacobian is
\[
J(i,x(i))= \begin{bmatrix} c_{11}(i) (1+x_1(i))^{-2} &c_{12}(i)(1+x_2(i))^{-2}\\
c_{21}(i)(1+x_1(i))^{-2} &c_{22}(i) (1+x_2(i))^{-2} \end{bmatrix} .
 \]
Calculation of the integral in Assumption~\ref{assum:trp} yields 
 \[
		 \int_0^ 1  J(i,r a+(1-r) b ) \diff r  = 
		\begin{bmatrix}
												c_{11}(i) ( 1+a_1+b_1+a_1 b_1 )^{-1} &  c_{12}(i) ( 1+a_2+b_2+a_2 b_2 )^{-1}\\
												c_{21}(i) ( 1+a_1+b_1+a_1 b_1 )^{-1} &  c_{22}(i) ( 1+a_2+b_2+a_2 b_2 )^{-1} 
		\end{bmatrix}.
\]

We assume that there exist~$\alpha,\beta,\gamma>0$ such that
\begin{align*}
\alpha< c_{pq}(i) <\beta \text{ and } c_{11}(i)c_{22}(i) -c_{12}(i)c_{21}(i) >\gamma
\end{align*}
for all~$1\leq p,q\leq 2$ and all~$i\geq 0$.
Then it is clear that~$\Omega:=[0,v]\times[0,v]$ is an 
invariant set of~\eqref{eq:exalp} for any~$v>0$ sufficiently large,
and that Assumption~\ref{assum:trp} holds. We also assume that the~$c_{pq}(i)$'s are all periodic with a common period~$T$. 
Then  Theorem~\ref{thm:nonltprt} implies that every solution converges to a~$T$-periodic trajectory. 

Figure~\ref{fig:imp}
depicts the solution of~\eqref{eq:exalp} for
\begin{align*}
c_{11}(i)&=5+\sin(\frac{i \pi }{2} +0.2),\\
c_{12}(i)&=2+\sin(\frac{i \pi }{2}   ),\\
c_{21}(i)&=3/2,\\
c_{22}(i)&=5+\sin( 2 i \pi    +4),
\end{align*}
and the initial condition~$x_0=\begin{bmatrix} 5&6 \end{bmatrix}'$.
 Note that for these values the  system is~$T$-periodic for~$T=4 $. It may be seen that the trajectory~$x(i)$
 indeed converges to a periodic pattern, with the same  period~$T$. 
\end{Example}

Consider now the  time-invariant nonlinear system
\be\label{eq:tins}
x(i+1)=f(x(i))
\ee
whose trajectories evolve on a compact and convex set~$\Omega\subset\R^n$. 
The associated variational system is
\begin{align}\label{eq:vstrew}
z(i+1)
&=\left ( \int_0^ 1  J( r x(i)+(1-r) \bar x(i)) \diff r \right ) 
 z(i).
\end{align}

The next assumption guarantees that~\eqref{eq:vstrew} is~TPDTS.
\begin{Assumption}\label{assum:trpnli}
						The matrix~$ \int_0^ 1  J( r a+(1-r) b ) \diff r $ is~TP for any~$a,b \in \Omega $.
\end{Assumption} 

Pick an arbitrary integer~$k\geq 0$. 
The time-invariant vector field is~$T$-periodic for~$T=k$ 
and thus entrainment implies that every solution converges to a~$k$-periodic trajectory. Since~$k$ is arbitrary this yields the following result. 
\begin{Corollary}\label{coro:eq}
If Assumption~\ref{assum:trpnli}  
 holds then any solution of~\eqref{eq:tins}
 converges to an equilibrium point. 
\end{Corollary}

This result may be regarded as the analogue of a    result of Smillie~\cite{smillie}
on convergence to an equilibrium
 in a continuous-time nonlinear cooperative systems with a special Jacobian. 

\section{Conclusions}
$SSR$ matrices appear in many fields. The most prominent 
example are~TP matrices.
Here we have studied firstly the spectral properties
 of matrices that are~$SSR_k$ for a single value~$k$. 
An important property of such matrices 
 is that some  eigenvalues  
can be complex (unlike in the $SSR$ case, where all eigenvalues are real).  
We then showed that the investigation of matrices that are~$SSR_k$ for a single value~$k$
can be used as
 a basic building block for studying matrices that are~$SSR_k$ for several values of~$k$.
e.g., $SSR$ matrices or matrices that are~$SSR_k$ for all odd~$k$.

 As an application, we derived an 
analogue  of the notion of~TPDS for the discrete-time case.
We introduced  the notion of a~TPDTS, that is,
 linear time-varying discrete-time system whose
transition matrices are all~TP. We showed that the vector solutions of such matrices have special sign properties. We then showed how this can be used to analyze the asymptotic behavior of nonlinear time-varying 
discrete-time systems whose
 variational equation is a~TPDTS.

 As explained in the recent paper~\cite{fulltppaper}, VDPs satisfied by 
the solutions of  linear time-varying
 systems can be used to prove the stability of certain nonlinear dynamical systems.
 In this context, it may be of interest to study the following problem. 
Consider the system
\be \label{eq:xdoo}
\dot x(s)=A(s)x(s),
\ee
where~$x(s)\in\R^n$ and~$A(s)$ is a continuous matrix function of~$s$. 
For any pair~$t_0\leq t$,
the solution of~\eqref{eq:xdoo} satisfies~$x(t)=\Phi(t,t_0)x(t_0)$, where~$\Phi(t,t_0)$
is the transition matrix from time~$t_0$ to time~$t$, that is, 
 the solution at time~$t$ of the matrix differential equation
\[
\dot \Phi(s)=A(s)\Phi(s), \text{ with }  \Phi(t_0)=I. 
\]
An interesting question is: given a value~$k \in\{1,\dots,n\}$,
when will~$\Phi(t,t_0)$ be~$SSR_k$ for all pairs~$t_0,t$ with~$t>t_0$, and what will 
be the implications of this for the solution of~\eqref{eq:xdoo}? 

It is well-known that the sum of two~TP matrices is not 
necessarily a~TP matrix. 
This means that establishing that Assumption~\ref{assum:trp}
indeed holds, that is, that 
						the matrix~$ \int_0^ 1  J(i,r a+(1-r) b ) \diff r $ is~TP for all~$a,b \in \Omega $ and all~$i\geq0$ is not trivial. An interesting direction for further
						research is to find simple conditions guaranteeing that this indeed holds. 

\bibliographystyle{IEEEtranS}
\bibliography{rola_bib}

 \end{document}